\documentclass[11pt]{amsart}
\usepackage{epsfig}
\usepackage{verbatim}
\usepackage{amssymb}
\newcommand{\ra}{{\rightarrow}}
\newcommand{\lra}{{\longrightarrow}}
\newcommand{\eproof}{\hfill\rule{2.2mm}{3.0mm}}

\newcommand{\Proof}{\noindent {\bf Proof.~~}}
\newcommand{\D}{{\mathcal D}}

\newcommand{\R}{{\mathbb R}}
\newcommand{\Z}{{\mathbb Z}}

\newcommand{\Q}{{\mathbb Q}}
\newcommand{\N}{{\mathbb N}}
\newcommand{\ve}{{\mathbf e}}
\newcommand{\vu}{{\mathbf u}}
\newcommand{\vv}{{\mathbf v}}
\newcommand{\ep}{\varepsilon}

\renewcommand{\eqref}[1]{(\ref{#1})}

\newcommand{\shsp}{\hspace{1em}}
\newcommand{\mhsp}{\hspace{2em}}
\newcommand{\FT}[1]{\widehat{#1}}

\newcommand{\diag}{{\rm diag}}
\newcommand{\J}{{\mathcal J}}
\renewcommand{\L}{{\mathcal L}}
\newcommand{\GLd}{{\rm GL}(d,\R)}
\newcommand{\GLtwo}{{\rm GL}(2,\R)}
\newtheorem{prop}{Proposition}[section]
\newtheorem{lem}[prop]{Lemma}

\newtheorem{theo}[prop]{Theorem}

\begin{document}
\baselineskip 18pt
\title{Simultaneous Translational and Multiplicative Tiling and Wavelet Sets in $\R^2$}
\author{Eugen J. Ionascu}
\address{Columbus State University\\4225 University Ave.\\
Columbus, GA 31907, USA, also Honorific Members of the Romanian
Institute of Mathematics ``Simion Stoilow''}
\email{ionascu@math.gatech.edu,ionascu$\_$eugen@colstate.edu}
\author{Yang Wang}
\thanks{Supported in part by the
       National Science Foundation grant DMS-0456538.}
\address{School of Mathematics  \\ Georgia Institute of Technology \\
Atlanta, GA 30332-0160, USA.}
\email{wang@math.gatech.edu}

\begin{abstract}
     Simultaneous tiling for several different translational sets has been studied
rather extensively, particularly in connection with the Steinhaus
problem. The study of orthonormal wavelets in recent years,
particularly for arbitrary dilation matrices, has led to the study
of multiplicative tilings by the powers of a matrix. In this paper
we consider the following simultaneous tiling problem: Given a
lattice in $\L\in \R^d$ and a matrix $A\in\GLd$, does there exist a
measurable set $T$ such that both $\{T+\alpha:~\alpha\in\L\}$ and
$\{A^nT:~n\in\Z\}$ are tilings of $\R^d$? This problem comes
directly from the study of wavelets and wavelet sets. Such a $T$ is
known to exist if $A$ is expanding. When $A$ is not expanding the
problem becomes much more subtle. Speegle \cite{Spe03} exhibited
examples in which such a $T$ exists for some $\L$ and nonexpanding
$A$ in $\R^2$. In this paper we give a complete solution to this
problem in $\R^2$.
\end{abstract}

\maketitle

\section{Introduction}
\setcounter{equation}{0}

The history of tiling goes as far back as the beginning of
civilization. Tiling has been studied in the history of mankind in
many different contexts and for different purposes.  Mathematically
we usually study tiling in the context of having a finite set of
``shapes''called {\em prototiles} and using congruent copies of
these prototiles to cover the whole Euclidean space without
overlapping. Translational tiling is one such example, in which only
the translations of the prototiles are used to tile the space.

Recently attentions have been given to {\em multiplicative tilings}. In a multiplicative
tiling there is a finite set of prototiles $\{T_1, T_2, \dots, T_m\}$ and sets
of nonsingular $d\times d$ matrices $\D_1, \D_2, \dots, \D_m$ such that
$$
    \Bigl\{A_jT_j:~A_j\in\D_j, 1\leq j \leq m\Bigr\}
$$
is a partition of $\R^d$. Here we define a partition in the most general sense, namely
the sets are disjoint in Lebesgue measure and their union is $\R^d$ up to a measure
zero set, see e.g. Wang \cite{Wan04} and Speegle \cite{Spe03}.
These studies are motivated in large part by the connection with orthonormal
wavelets, which we shall discuss later.

We first introduce some notations and terminologies. We say a collection of measurable
sets $\{T_j\}$ in $\R^d$ is a {\em tiling} of
$\R^d$ if it is a partition of $\R^d$ in the sense
just described above. A measurable set $T$ is said to tile {\em translationally by
$\J$}, where $\J$ is a subset of $\R^d$, if $\{T+\alpha:~\alpha\in\J\}$
is a tiling of $\R^d$. In this paper, we are primarily concerned with translational
tiling by a lattice. Let $A\in \GLd$ be a $d\times d$ nonsingular matrix. We say
a measurable set $T$ {\em tiles multiplicatively by $A$}
if $\{A^n\,T:~n\in\Z\}$ is a tiling of $\R^d$. The main question we ask in this paper
is:

\vspace{2mm}
\noindent
{\bf Problem.}~~{\em Given a matrix $A\in\GLd$ and a lattice $\L$ in
$\R^d$, does there exist a measurable set $T$ in
$\R^d$ such that $T$ tiles translationally by $\L$ and multiplicatively
by $A$?}

Simultaneous translational tiling using more than one lattice have been studied rather
extensively. One of the best known problem in tiling is the classic
Steinhaus problem posed by H. Steinhaus sometime in the 1950's, which asks
for the existence of a $T\subset\R^2$ that tiles translationally by all lattices of the
form $R_\theta\Z^2$ where $\R_\theta$ is a rotation matrix. It was shown by
Jackson and Mauldin \cite{JaMa02} that such a $T$ exists, but in their construction
$T$ is not measurable. The problem remains open for measurable sets. Han and Wang
\cite{HaWa01} proved that for any two lattices $\L_1,\L_2$ in $\R^d$ with the same
co-volume there always exists a measurable $T\subset\R^d$ such that $T$ tiles translationally
by both $\L_1$ and $\L_2$. This problem is motivated by the study of Weyl-Heisenberg
orthonormal bases for $L^2(\R^d)$. Kolountzakis \cite{Kol97} established conditions on lattices
$\L_1, \L_2, \dots, \L_n$ in $\R^d$ with the same co-volume for which a measurable
$T$ exists that tiles translationally by each $\L_j$. There are many other related studies
on translational simultaneous tilings, see the references in the aforementioned papers.

In contrast, the study of simultaneous translational and multiplicative tiling has only
begun very recently.
One of the motivations for studying this problem is the connection to orthonormal wavelets.
For $A\in\GLd$ we call a function $f\in L^2(\R^d)$ an {\em orthogonal wavelet with dilation
$A$} if the set of functions
$$
    \Bigl\{|\det(A)|^{\frac{n}{2}}f(A^n x-\alpha):~n\in\Z, \alpha\in\Z^d\Bigr\}
$$
is an orthogonal basis for $L^2(\R^d)$. Again, in general we can
substitute $\Z^d$ with any full rank lattice. Since the seminal work
of Daubechies \cite{Dau88}, Lemari$\acute{e}$ and Mayer \cite{LM86}
there has been an explosion in the study of wavelets and their
applications in image compression, digital signal processing and
numerical computations. We shall not go into details about wavelets
here and shall refer the readers to Daubechies \cite{Dau-book} and
Mayer \cite{M92}. It should be pointed out that in most studies on
wavelets the dilation matrix $A$ is assumed to be {\em expanding},
i.e. all eigenvalues have $|\lambda|>1$. Furthermore, a single
dilation matrix is involved. In \cite{Wan04}, the concept of
wavelets is broadened to allow more than one dilation matrix as well
as nonexpanding matrices.

The role of tiling has appeared in the study of the following fundamental question in wavelets:
Given an expanding matrix $A\in\GLd$ is it always possible to find an orthogonal
wavelet with dilation $A$? This question was answered by studying functions of
the form $f = \FT \chi_T$ where $T$ is a subset of $\R^d$ with finite measure. Dai and
Larson \cite{DaLa98}\footnote{In their theorem the matrix $A$ is assumed to be expanding,
but it is clearly not needed.}proved the following theorem:
\begin{theo}[\cite{DaLa98}]  \label{theo-1.1}
Let $A\in\GLd$ and let $T\subset\R^d$ have finite measure. Then
$f=\FT \chi_T$ is an orthogonal wavelet with dilation $A$ if and
only if $T$ tiles $\R^d$ translationally by $\Z^d$ and
multiplicatively by $A^T$.
\end{theo}
A set $T$ that tiles $\R^d$ translationally by $\Z^d$ and
multiplicatively by $A^T$ simultaneously is called a {\em wavelet
set with dilation $A$}. Later, it was proved by Speegle, Dai and
Larson \cite{DLS97} that for an expanding or contracting $A$ a
wavelet set that is bounded exists. Many other different
constructions in the expanding case have been proposed, all of
which involve cut and paste, see e.g.
\cite{BeLe99,DLS97,ILP98,SoWe98,Wan04}. One starts with a set
$T_0$ that tiles by $A^T$ and covers $\R^d$ translationally by
$\Z^d$. The goal then is to move pieces of $T_0$ around using a
combination of translations and dilations to get a wavelet set.
The fact that $A$ is expanding or contracting plays a crucial role
because it can be used multiplicatively to control the size of the
pieces, an important part of the constructions.

The existence of wavelet sets for nonexpanding matrices is a much
more challenging problem. Wang \cite{Wan04} exhibited an example of
a wavelet set whose dilations consist more than the powers of a
single matrix, in which some dilations are neither expanding nor
contracting. But few people believed, or even thought about, the
possibility of an orthogonal wavelet with a single dilation matrix
$A$ that is neither expanding or contracting. However, Speegle
\cite{Spe03} recently has shown the existence of such wavelet sets
(and hence orthogonal wavelets). In particular, he has shown that
for $A={\rm diag} \,(a, b)$ where $|b|<1$ and $|ab|>1$ and the
lattice $\L=\Z \ve_1+\Z (\ve_2+\sqrt 3 \ve_1)$ there exists a
measurable $T$ that tiles $\R^d$ translationally by $\Z^2$ and
multiplicatively by $A$. More generally, the value $\sqrt 3$ can be
replaced with any $\beta \in\R$ if $\beta$ cannot be approximated by
rationals to within certain order $J=J(a,b)$. By a simple linear
transformation, one can show the existence of a matrix $A\in \GLd$
for $d=2$ which is neither expanding nor contracting, for which
there is a wavelet set.

In this paper we give a complete classification in dimension $d=2$
for the existence of wavelet sets, or
equivalently the existence of simultaneous tiling
translationally by a lattice and multiplicatively by $A$.
We state our main theorems here.

\begin{theo}  \label{theo-1.2}
Let $A\in \GLtwo$ with $|\det(A)| \geq 1$. Let $\lambda_1,\lambda_2$ be the eigenvalues
of $A$ with $|\lambda_1| \geq |\lambda_2|$.
\begin{itemize}
\item[\rm (a)] If $|\lambda_1\lambda_2|=1$, i.e. $|\det(A)|=1$, then there
exists no wavelet set with dilation $A$.
\item[\rm (b)] If $|\lambda_1| >1$ and $|\lambda_2|\geq 1$, then there
exists a wavelet set with dilation $A$.
\item[\rm (c)] If $|\lambda_1\lambda_2| >1$ and $|\lambda_2|< 1$, let $\vv=[v_1, v_2]^T$ be
an eigenvector of $A^T$ for $\lambda_2$. Then there exists a wavelet set
with dilation $A$ if and only if $v_1/v_2 \in\R\setminus\Q$.
\end{itemize}
\end{theo}

Note that the assumption $|\det(A)| \geq 1$ is without any loss of generality. If a
wavelet set exists with dilation $A$ then it is also a wavelet set with dilation $A^{-1}$.
We also remark that part (a) of the theorem is due to \cite{LSST05}, where it
is shown that if $|\det(A)|=1$ then there exists no $T$ with finite measure such that
$T$ tiles multiplicatively by $A$.
Furthermore, they proved that a {\em bounded} wavelet set exists if and only if $A$ is
expanding or contracting. Our next theorem is more general than Theorem \ref{theo-1.2},
and it gives a complete classification of the simultaneous tiling problem in $\R^2$.

\begin{theo}  \label{theo-1.3}
Let $A\in \GLtwo$ with $|\det(A)| \geq 1$ and eigenvalues
$\lambda_1,\lambda_2$, $|\lambda_1| \geq |\lambda_2|$. Let
$\L=P\Z^2$ be a lattice in $\R^2$ with $P\in\GLtwo$.
\begin{itemize}
\item[\rm (a)] If $|\lambda_1\lambda_2|=1$, then there exists no measurable $T \subset\R^2$
such that $T$ tiles translationally by $\L$ and multiplicatively by $A$.
\item[\rm (b)] If $|\lambda_1| >1$ and $|\lambda_2|\geq 1$, then there
exists a measurable $T \subset\R^2$
such that $T$ tiles translationally by $\L$ and multiplicatively by $A$.
\item[\rm (c)] If $|\lambda_1\lambda_2| >1$ and $|\lambda_2|< 1$, let $\vv=[v_1, v_2]^T$ be
an eigenvector of $P^{-1}AP$ for $\lambda_2$. Then there exists a
measurable $T \subset\R^2$ such that $T$ tiles translationally by
$\L$ and multiplicatively by $A$ if and only if $v_1/v_2
\in\R\setminus\Q$.
\end{itemize}
\end{theo}

This work is completed while the first author is visiting the School of Mathematics
of Georgia Institute of Technology as part of the Faculty Development Program
sponsored by the University System of Georgia. He would like to thank both the
Columbus State University and Georgia Tech for the generous support.

\section{Some General Results}
\setcounter{equation}{0}

We focus here on general results concerning simultaneous multiplicative and translational
tilings. These results will be used to prove our main theorems. First we introduce some
terminologies.

Let $\Omega\subset\R^d$ be a measurable set. We say {\em $\Omega$ packs $\R^d$ translationally
by $\L$} or simply {\em $\Omega$ packs by $\L$}, where
$\L$ is a lattice in $\R^d$, if $\{\Omega+\alpha:~\alpha\in\L\}$ are disjoint in measure.
Similarly, we say
{\em $\Omega$ packs $\R^d$ multiplicatively
by $A$} or simply {\em $\Omega$ packs by $A$}, where
$A\in\GLd$, if $\{A^n\Omega:~n\in\Z\}$ are disjoint in measure.
Since the construction of wavelet sets and tiles often involve cut and paste, the concept of
packing plays an important role in this paper.

\begin{lem}  \label{lem-2.1}
   Let $A\in\GLd$ and $\L = P\Z^d$ with $P\in\GLd$. Then a set $T$ tiles multiplicatively by $A$
   and translationally by $\L$ if and only if $P^{-1}T$ tiles multiplicatively by $P^{-1}AP$
   and translationally by $\Z^d$.
\end{lem}
\Proof   If $\{A^nT:~n\in\Z\}$ is a tiling then clearly so is
$\{P^{-1}A^nT=(P^{-1}AP)^nP^{-1}T:~n\in\Z\}$, and conversely. Similarly, if
$\{T+P\alpha:~\alpha\in\Z^d\}$ is a tiling then so is
$\{P^{-1}T+\alpha:~\alpha\in\Z^d\}$, and conversely.
\eproof

In the study of multiplicative tiling it is useful to study the
{address map}. Let $\Omega$ be a multiplicative tile by $A\in\GLd$.
Then the {\em address map (induced by $\Omega$)} is the map
$\tau_\Omega: \R^d \lra \Omega$ given by
\begin{equation}  \label{2.1}
    \tau_\Omega(x) = y  \mhsp\mbox{if}\mhsp x=A^n y, y\in\Omega.
\end{equation}
Note that for almost all $x\in\R^d$ there exist
unique $n\in\Z$ and $y\in\Omega$ such
that $x=A^ny$, so the map is well defined almost everywhere. If $S\subset\R^d$ packs
by $A$ then $\tau_\Omega|_S$ is one-to-one almost
everywhere on $S$, and conversely. Furthermore,
for such an $S$ we can define the {\em index map (induced by $S$ and $\Omega$)}
$\phi:~\Omega \lra \Z$ given by
\begin{equation}  \label{2.2}
    \phi(y) =\left\{\begin{array}{ll} n & \shsp\mbox{if $A^n y\in S$ for some $n\in\Z$} \\
                                      0 & \shsp\mbox{otherwise.}
            \end{array}\right.
\end{equation}

\begin{theo}  \label{theo-2.2}
    Let $A\in\GLd$ and $\L$ be a lattice in $\R^d$. Suppose there exists an $S\subset\R^d$
    such that $S$ tiles multiplicatively by $A$ and packs translationally by $\L$.
    Then there exists a $T\subset\R^d$
    such that $T$ tiles multiplicatively by $A$ and translationally by $\L$.
\end{theo}
\Proof  Note that we must have $|\det(A)| \neq 1$, for otherwise $S$
would have infinite measure and cannot pack translationally by $\L$,
see Theorem~4 in \cite{LSST05}. It is also proved in \cite{LSST05}
that if $|\det(A)| \neq 1$ then we can find an $S_0\subset\R^d$ such
that $S_0$ tiles multiplicatively by $A$ and furthermore, the
construction in the paper clearly shows that we can require $S_0$ to
have nonempty interior. Now, set $S^*=nS_0$ for some $n$
sufficiently large. Again $S^*$ tiles multiplicatively by $A$. But
now $S^*$ contains a set $F$ that tiles by $\L$. Also, since $S$
packs by $\L$, there exists an $F^*$ that tiles by $\L$ with
$S\subseteq F^*$.

Now since both $F, F^*$ tiles translationally by $\L$, there exists
a bijection (in the sense of almost everywhere) $\rho: F^* \lra F$
such that $\rho(x)=x+\alpha(x)$ for some unique $\alpha(x)\in\L$.
The map $\tau_{S}: S^* \lra S$ is also a bijection since both $S$,
$S^*$ tile multiplicatively by $A$. Let $\phi= \rho|_S$ and
$\psi=\tau_S|_F$. Then both $\phi: S\lra F$ and $\psi: F\lra S$ are
one-to-one.

By the Schr\"oder-Cantor-Bernstein construction there exists a measurable bijection
$h: S \lra F$ having the form
$$
    h(x) =\left\{\begin{array}{ll} \phi(x) &\shsp x\in E \\
                        \psi^{-1}(x) &\shsp x\in S\setminus E \end{array}\right.
$$
for some $E\subseteq S$. A more precise form of $E$ can be found in \cite{Ion02}.
Clearly for each $x\in S$ there exist a unique $n(x)\in \Z$ and a unique
$\alpha(x)\in\L$
such that $h(x) = A^{n(x)}x+\alpha(x)$. Let $T=\{A^{n(x)}x: ~x\in S\}$. It is obvious
that $T$ tiles multiplicatively by $A$ because $S$ does. Furthermore, $T$ is congruent to
$F$ modulo the lattice $\L$, so $T$ also tiles by $\L$.
\eproof

The above theorem is established in special forms in \cite{IoPe98} and \cite{Spe03}.
It shows that to prove the existence of simultaneous multiplicative and translational
tilings one only needs to prove the existence of simultaneous multiplicative tiling and
translational packing. This is precisely the strategy we follow to prove our main theorems.

\begin{lem}  \label{lem-2.3}
Let $\Omega\subset\R^d$ with $\mu(\Omega) < \infty$ tile by $A \in\GLd$,
$|\det(A)|>1$. Let $S_n$ pack $\R^d$ by $A$ and
$\phi_n: \Omega \lra\Z$ be the index map induced by $S_n$ and $\Omega$.
Assume that $\lim_{n\ra\infty} \mu(S_n)=0$ and
$\liminf_n \phi_n(x) > -\infty$ for almost all $x\in\Omega$.
Then $\lim_{n\ra\infty}\mu(\tau_\Omega(S_n))=0$.
\end{lem}
\Proof  Assume that $\mu(\Omega)$ and $\mu(S_n)$ are all bounded by
$c>0$. Each $S_n$ has a unique decomposition $S_n = \bigcup_{m\in\Z}
A^m R_{n,m}$ (up to a null set) where $R_{n,m}=\Omega \cap A^{-m}S_n
\subseteq \Omega$. Furthermore, $\tau_\Omega$ is injective on $S_n$
so the sets $\{R_{n,m}\}$ are disjoint in measure, and
$\tau_{\Omega}(S_n) =\bigcup_{m\in\Z} R_{n,m}$. We have
\begin{equation} \label{2.3}
    \sum_{m\in\Z} \mu(R_{n,m}) \leq \mu(\Omega) \leq c, \mhsp
    \sum_{m\in\Z} \beta^m \mu(R_{n,m}) = \mu(S_n) \leq c,
\end{equation}
where $\beta :=|\det(A)|>1$.
For any $\ep>0$, note that $\sum_{m\geq N}\beta^m \mu(R_{n,m}) \leq c$ so
$$
      \sum_{m\geq N} \mu(R_{n,m}) \leq c\beta^{-N}.
$$
Hence for all $N \geq N_0(\ep)$ we have $\sum_{m\geq N} \mu(R_{n,m}) < \ep$.

Since $\liminf_n \phi_n(x) > -\infty$ for almost all $x\in\Omega$,
there exists an $N_1=N_1(\ep)$
such that for all $N, n \geq N_1$ we have
$$
    \mu(\{x\in\Omega:~\phi_n(x) <-N\}) < \ep,
$$
which is equivalent to $\sum_{m <-N} \mu(R_{n,m}) < \ep$. Now pick $N \geq N_0, N_1$
and we have
$$
    \sum_{|m|> N} \mu(R_{n,m})= \sum_{m >N} \mu(R_{n,m}) + \sum_{m <-N} \mu(R_{n,m})
       <2\ep.
$$
In addition,
$$
   \sum_{|m| \leq N} \beta^m \mu(R_{n,m}) \leq \mu(S_n)\lra 0
                                  \shsp \mbox{as} \shsp n \ra \infty.
$$
So there exists an $N_2$ such that $\sum_{|m| \leq N} \mu(R_{n,m})< \ep$
for all $n \geq N_2$. Thus for $n \geq N_2$ we have
$$
  \sum_{m\in\Z} \mu(R_{n,m}) <3\ep.
$$
It follows that $\mu(\tau_{\Omega}(S_n)) = \sum_{m\in\Z} \mu(R_{n,m}) \lra 0$
as $n \ra \infty$.
\eproof

\begin{lem}  \label{lem-2.4}
Let $\Omega\subset\R^d$ with $\mu(\Omega) < \infty$ tile by $A \in\GLd$,
$|\det(A)|>1$. Let $S_n$ pack $\R^d$ by $A$ and
$\phi_n:\Omega \lra\Z$ be the index map induced by $S_n$ and $\Omega$.
Assume $ S^* \subset\R^d$ with $\mu( S^*)<\infty$ such that
$\lim_{n\ra\infty} \mu(S_n\Delta  S^*)=0$, where $\Delta$ denotes the symmetrical
difference. Then
\begin{itemize}
\item[\rm (a)]  $ S^*$ packs $\R^d$ multiplicatively by $A$.
\item[\rm (b)]  Assume that $\lim_{n\ra\infty} \mu(\tau_\Omega(S_n)\Delta \Omega)=0$
and $\liminf_n \phi_n(x) > -\infty$ for almost all $x\in\Omega$.
Then $ S^*$ tiles $\R^d$ multiplicatively by $A$.
\end{itemize}
\end{lem}
\Proof  (a) is stated in \cite{Spe03} without a proof. It is quite
straightforward, but we will furnish a proof here. Assume it is
false, then $\mu(A^k S^*\cap A^l S^*) = \delta>0$ for some $k \neq
l$. But $\lim_n\mu(S_n \Delta S^*) =0$. It follows that
$\mu(A^kS_n\cap A^lS_n) \geq \delta/2>0$ for sufficiently large $n$,
a contradiction.

We now prove (b) by proving $\tau_\Omega( S^*)=\Omega$. Assume it is
false. Then $\mu(\Omega\setminus\tau_\Omega( S^*))>0$. Hence, by the
first hypothesis in (b), there exists a $\delta>0$ such that
$\mu(\tau_\Omega(S_n)\setminus\tau_\Omega( S^*))>\delta$ for all
sufficiently large $n$. Thus $\mu(\tau_\Omega(S_n\setminus S^*)
>\delta$ for sufficiently large $n$. We shall derive a
contradiction.

Set $R_n = S_n\setminus S^*$. Then $\lim_n \mu(R_n)=0$. Note that by
assumption $\liminf_n \phi_n(x) > -\infty$ for almost all $x\in\Omega$. Let
$\psi_n$ be the index map induced by $R_n$ and $\Omega$. Then either $\psi_n(x)=
\phi_n(x)$ or $\psi_n(x) = 0$ because $R_n \subseteq S_n$. Thus
$\liminf_n \psi_n(x) > -\infty$ for almost all $x\in\Omega$. It follows from
Lemma \ref{lem-2.3} that $\lim_n \mu(\tau_\Omega(R_n)) = 0$, a contradiction.
\eproof

The next theorem is a stronger and more general
version of Theorem 3.2 in \cite{Spe03}. We also
give a different proof here, using the above lemmas.

\begin{theo}  \label{theo-2.5}
Let $\L$ be a lattice in $\R^d$ and $A\in\GLd$ with $|\det(A)|>1$. Suppose that
for any bounded set $S\subset\R^d$ there are infinite many $n\in\N$ such that
$A^{-n}S$ packs translationally by $\L$. Then there exists
a $T\subset \R^d$ such that $T$ tiles $\R^d$ multiplicatively by $A$ and
translationally by $\L$.
\end{theo}
\Proof By \cite{LSST05} there exists an $\Omega$ with $\mu(\Omega)=c<\infty$ and
$\Omega$ tiles multiplicatively by $A$. If $\Omega$ is also bounded then we can find
$m_1>0$ such that $A^{-m_1}\Omega$ packs by $\L$. Since $A^{-m_1}\Omega$ also
tiles by $A$, the theorem follows from Theorem \ref{theo-2.2}. Thus we shall
assume $\Omega$ is unbounded.

Denote $\Omega_n =\Omega\cap [-M_n,M_n]^d$, where
$0<M_1<M_2<\cdots$ with the property that $\mu(\Omega\setminus
\Omega_n) \leq 4^{-n}c$. Write $\Omega_0= \emptyset$ and $R_n =\Omega_n\setminus
\Omega_{n-1}$ for $n \geq 1$. Then $\{R_n\}$ are disjoint and
$$
    \Omega = \bigcup_{k=1}^\infty R_k, \mhsp
    \Omega_n = \bigcup_{k=1}^n R_k.
$$
The idea is to use $\{R_n\}$ to construct a sequence of sets $\{S_n\}$ satisfying
the conditions of Lemma \ref{lem-2.4} (b).

Let $S_1 = A^{-m_1}\Omega_1$ such that $S_1$ packs by $\L$. This $m_1$ exists since $\Omega_1$
is bounded. We construct recursively $S_n$ for $n \geq 2$
satisfying the following properties:
\begin{itemize}
\item[(a)]  $S_n$ packs by $\L$ and $\tau_\Omega(S_n) = \Omega_n$.
\item[(b)]  $\mu(S_n\Delta S_{n-1}) \leq 4^{-n}\mu(S_1)$.
\item[(c)]  Each $\phi_n$ is bounded and the set
$X_n:=\{x\in\Omega:~\phi_n(x) \neq \phi_{n-1}(x)\}$ has
$\mu(X_n) \leq 2c/4^{n}$, where $c=\mu(\Omega)$ and
$\phi_n$ is the index map induced by $S_n$ and $\Omega$.
\end{itemize}

Assume such  $S_n$'s exist. Observe that for all $n>m$ we have
\begin{equation}  \label{2.4}
\mu(S_n\Delta S_m) \leq \sum_{k=m}^{n-1} \mu(S_{k+1}\Delta S_k)
                   \leq \sum_{k=m}^{n-1} \frac{1}{4^{k+1}}\mu(S_1) <\frac{2}{4^m}\mu(S_1).
\end{equation}
So $\{S_n\}$ is a Cauchy sequence in the sense of symmetrical difference, and there
exists an $ S^*$ such that $\lim_n \mu(S_n\Delta S^*) = 0$.
Furthermore, taking $m=1$ in (\ref{2.4}) yields
$\mu(S_n\Delta S_1) < \mu(S_1)/2$, so $\mu(S_n) > \mu(S_1)/2$. This means
$ S^*$ has positive measure. By Lemma 3.1 in \cite{Spe03}, $S^*$ packs
by $A$ and $\L$. Now let
$$
    Y_n:=\Bigl\{x\in\Omega:~\phi_n(x) \neq \phi_k(x) \mbox{~for all~}k>n\Bigr\}.
$$
Then (c) yields
$$
 \mu(Y_n) \leq \sum_{k>n} \mu(X_k) \leq \sum_{k>n} \frac{2c}{4^k} <\frac{c}{4^n}.
$$
Since each $\phi_n$ is bounded on $\Omega$, it follows that
$\liminf_n \phi_n(x) >-\infty$ for almost all $x\in\Omega$. Hence $ S^*$ tiles multiplicatively
by $A$. The theorem follows from Theorem \ref{theo-2.2}.

It remains to prove that such $S_n$'s exist. Assume $S_{n-1}$ has been constructed,
$n>1$. We construct $S_n$. Let $m_n>0$ such that $A^{-m_n}\Omega_n$ packs by $\L$.
Let
$$
     U_n = S_{n-1} \cap \bigl(A^{-m_n}\Omega_n +\L\bigr).
$$
Since both $S_{n-1}$ and $A^{-m_n}\Omega_n$ pack by $\L$, it is easy to see that
$$
    \mu(U_n) \leq \mu(A^{-m_n}\Omega_n) \leq |\det(A)|^{-m_n}\, c.
$$
Now $\phi_{n-1}$ is bounded, so $\phi_{n-1}(x) \geq -a_{n-1}$ for some $a_{n-1}>0$.
This means that for each $x\in S_{n-1}$ there exists a $k\leq a_{n-1}$ such that
$x = A^{-k}\tau_\Omega(x)$. Hence
$\mu(\tau_\Omega(U_n)) \leq |\det(A)|^{a_{n-1}}\mu(U_n)$.
By choosing $m_n \geq a_{n-1}$
sufficiently large we have $\mu(\tau_\Omega(U_n)) \leq 4^{-n}c$ and
$\mu(A^{-m_n}\Omega_n) \leq 4^{-n}\mu(S_1)$. We now define
$$
     S_n =\left( S_{n-1}\setminus U_n\right) \cup
          A^{-m_n}\tau_\Omega(U_n) \cup A^{-m_n}R_n.
$$
It is clear that $S_n$ packs by $\L$, and $\tau_\Omega(S_n) =\tau_\Omega(S_{n-1})\cup R_n
=\Omega_n$. So (a) is satisfied. Also,
$$
     \mu(S_n\Delta S_{n-1}) \leq \mu(U_n) + \mu(A^{-m_n}(\tau_\Omega(U_n) \cup R_n))
                            \leq 2\mu(A^{-m_n}\Omega_n)
                            \leq \frac{\mu(S_1)}{4^{n}}.
$$
So (b) is satisfied. Finally, observe that $\phi_n(x) = -m_n$ for $x\in\tau_\Omega(U_n) \cup R_n$
but $\phi_n(x)=\phi_{n-1}(x)$ everywhere else. So $\phi_n$ is bounded. Furthermore,
$$
    \mu(X_n) \leq \mu(\tau_\Omega(U_n) \cup R_n) <\frac{2c}{4^n}.
$$
This yields (c). The proof of the theorem is now complete.
\eproof

\section{Proof of Main Theorems}
\setcounter{equation}{0}

   We now prove our main theorems. The proofs are divided into several propositions for
different cases. A key ingredient is the approximation of irrational numbers
by rational numbers.

We first consider the case in which $A\in\GLtwo$ has eigenvalues $\lambda_1, \lambda_2$
with $|\lambda_1| >|\lambda_2|=1$. In this case, both $\lambda_1,\lambda_2$
are necessarily real, and $\lambda_2 = \pm 1$.

\begin{prop}   \label{prop-3.1}
     Let $A\in\GLtwo$ with eigenvalues $\lambda_1, \lambda_2$,
$|\lambda_1| >|\lambda_2|=1$. Assume $A$ has no rational eigenvectors
for $\lambda_2$. Then there exists a measurable $T\subset\R^2$ such that
$T$ tiles translationally by $\Z^2$ and multiplicatively by $A$.
\end{prop}
\Proof  Let $\vv_1$ and $\vv_2$ be eigenvectors of $A$ for $\lambda_1$ and
$\lambda_2$, respectively. By assumption we may assume $\vv_2=[1, \beta]^T$
where $\beta\in\R\setminus\Q$. We use Theorem \ref{theo-2.5} to complete the proof.
Let $S$ be any bounded set in $\R^2$. We prove $A^{-n}S$ packs translationally
by $\Z^2$ for all sufficiently large $n$. Since $S\subseteq
\{s\vv_1+t\vv_2:~|s|,|t|\leq K\}$ for some $K>0$, we may without loss of generality
assume $S=\{s\vv_1+t\vv_2:~|s|,|t|\leq K\}$. Observe that
$$
    A^{-n}S=\{s\vv_1+t\vv_2:~|s|\leq |\lambda_1|^{-n}K,|t|\leq K\},
$$
which approaches the straight line segment $L = \{t\vv_2:~|t| \leq
K\}$ in Hausdorff metric.

Now $\vv_2=[1,\beta]^T$ has irrational slope,
so the line segments $\{L+\alpha:~\alpha\in\Z^2\}$
are disjoint. Let $\ep_0 = {\rm dist}\,(L, \Z^2\setminus\{0\})>0$.
Then the distance between $L+\alpha_1$ and $L+\alpha_2$ is at least $\ep_0$ for any
$\alpha_1 \neq \alpha_2$. It follows that for sufficiently large $n$, by making
$|\lambda_1|^{-n}K< \ep_0/3$, the sets
$\{ A^{-n}S+\alpha:~\alpha\in\Z^2\}$ are disjoint. Thus $A^{-n}S$ packs by $\Z^2$ for
sufficiently large $n$, proving the proposition.
\eproof

\begin{prop}   \label{prop-3.2}
     Let $A\in\GLtwo$ with eigenvalues $\lambda_1, \lambda_2$,
$|\lambda_1| >|\lambda_2|=1$. Assume $A$ has a rational eigenvector
for $\lambda_2$. Then there exists a measurable $T\subset\R^2$ such that
$T$ tiles translationally by $\Z^2$ and multiplicatively by $A$.
\end{prop}
\Proof  Let $\vv_2=[p,q]^T$ be an eigenvector of $A$ for $\lambda_2$,
with $p,q\in\Z$ and $\gcd\, (p,q)=1$. Let $\vv_1 \in\Z^2$ with the property that
$P=[\vv_1,\vv_2]$ has $\det(P)=1$, where the columns of $P$ are $\vv_1$ and $\vv_2$.
For this $P$ we have
$$
P^{-1}AP = \begin{bmatrix} \lambda_1 & 0\\b & \lambda_2\end{bmatrix},\mhsp
    P\Z^2 = \Z^2.
$$
Now, take $Q=\bigl[\begin{smallmatrix} 1&0\\t&1\end{smallmatrix}\bigr]$ with
$t = \frac{b}{\lambda_1-\lambda_2}$. Then
$Q^{-1}P^{-1}APQ = \diag(\lambda_1, \lambda_2)$. Let $B=\diag(\lambda_1, \lambda_2)$.
By Lemma \ref{lem-2.1} it suffices to prove there exists a $T$ such that
$T$ tiles translationally by $Q^{-1}P^{-1}\Z^2=Q^{-1}\Z^2$
and multiplicatively by $B$.

We prove the existence by an explicit construction. Denote $U_n=[-a_n, a_{n+1})
\cup (a_{n+1},a_n]$ where $a_n =\frac{1}{2}|\lambda_1|^{-n}$. It is clear
that $\lambda_1^{-1} U_n = U_{n+1}$, $\{U_n\}$ is a partition of $\R$, and
$\bigcup_{n=0}^\infty U_n =[-1/2, 1/2]$ up to a null set. Let $I_n =
[-\frac{1}{2}-n, -n] \,\cup\, [n, n+\frac{1}{2}]$ and
$ \Omega_n= U_n \times I_n$, $n \geq 0$.
Observe that each $I_n$ tiles $\R$ translationally by $\Z$. So
$\Omega_n +\{0\}\times\Z =U_n \times \R$. Set $T= \bigcup_{n=0}^\infty \Omega_n$.

Now, $B^{-k}\Omega_n =U_{n+k}\times I_n$. Thus $\bigcup_{k\in\Z} B^{-k}\Omega_n
= \R\times I_n$ with the union disjoint. So $\bigcup_{k\in\Z} B^{-k} T=\R^2$
with the union disjoint. Hence $T$ tiles by $B$. It remains to prove $T$ tiles
$\R^2$ by $Q^{-1}\Z^2$. This follows from the observation that
$\{T+ [0,k]^T:~k\in\Z\}$ is a partition of $[-1/2, 1/2]\times\R$. Hence
$T$ tiles translationally by $Q^{-1}\Z^2$.
\eproof

The more difficult part of our theorems concerns the case in which $|\lambda_1|>1$
and $|\lambda_2|<1$. In this case again both $\lambda_1, \lambda_2\in\R$, so there
exists a $P\in\GLtwo$ such that $P^{-1}AP = \diag\,(\lambda_1, \lambda_2)$. It thus suffices
to consider a diagonal matrix and ask which lattices $\L$ lead to simultaneous translational
and multiplicative tilings by $\L$ and $A$.

\begin{lem}  \label{lem-3.3}
    Let $\beta\in\R\setminus\Q$ and let $p_n/q_n$ be the $n$-th convergent of
the continued fraction expansion of $\beta$, $q_n>0$. Denote $M_n
=q_n-1$. Let $c, \ep>0$. Then there exists an $n_0=n_0(c, \ep)$ such
that if $n \geq n_0$ then for any $p,q\in\Z$ with $\gcd(p,q)=1$, $1
\leq q \leq M_n$ we have
$$
   \left|\beta -\frac{p}{q}\right| \geq \frac{c}{qM_n^{1+\ep}}.
$$
\end{lem}
\Proof  Choose $n_0$ so that for all $n\geq n_0$ we have
$c/M_n^{\ep} <1/2$, and $q_k \leq q_n-2$ for all $1\leq k <n$. Assume the lemma
is false, then there exist $p^*,q^*\in \Z$, $\gcd(p^*,q^*)=1$ and $1 \leq q^* \leq M_n$
such that
\begin{equation}  \label{3.1}
   \left|\beta -\frac{p^*}{q^*}\right| < \frac{c}{q^*M_n^{1+\ep}}
   < \frac{1}{2q^*M_n} \leq \frac{1}{2{q^*}^2}.
\end{equation}
It follows that $p^*/q^* = p_k/q_k$ for some $k<n$, see e.g. \cite{Khi-book}. But we know
from the properties of continued fractions that
$$
   \left|\beta -\frac{p^*}{q^*}\right| =\left|\beta -\frac{p_k}{q_k}\right|
   \geq \frac{1}{q_k(q_{k+1}+q_k)},
$$
see also \cite{Khi-book}.
Now $q_k \leq q_n-2$, so $q_{k+1}+q_k \leq 2M_n$. Thus
$ \left|\beta -\frac{p^*}{q^*}\right| \geq \frac{1}{2q^*M_n}$, which contradicts
(\ref{3.1}).
\eproof

\begin{prop}   \label{prop-3.4}
    Let $A=\diag\,(\lambda_1, \lambda_2)$ with $|\lambda_1|>1>|\lambda_2|$ and
$|\lambda_1\lambda_2|>1$. Let $\L=\Z\vu+\Z\vv$ be a full rank lattice in $\R^2$
with $\vu=[u_1,u_2]^T$ and $\vv=[v_1,v_2]^T$. Assume that $u_1/v_1\in\R\setminus\Q$.
Then there exists a measurable $T\subset\R^2$ such that
$T$ tiles translationally by $\L$ and multiplicatively by $A$.
\end{prop}
\Proof  Write $a=|\lambda_1|$ and $b=|\lambda_2|^{-1}$. Then we
have $a>b>1$. We use Theorem \ref{theo-2.5} and prove that for any
bounded $S\subset\R^2$ there exist infinitely many $n>0$ such that
$A^{-n}S$ packs by $\L$. Assume this is false then there exists a
bounded $S\subset\R^2$ such that $A^{-n}S$ does not pack by $\L$
for all sufficiently large $n$. We derive a contradiction.

Without loss of generality we assume that $S = [-r, r]^2$ for some
$r$. Hence $A^{-n}S=[-ra^{-n},ra^{-n}]\times [-rb^n,rb^n]$. Assume
that $A^{-n}S$ does not pack translationally by $\L$. Since $\L$
is a lattice, we can find an $\alpha\neq 0$ in $\L$ such that
$A^{-n}S\cap (A^{-n}S+\alpha) \neq \emptyset$. Therefore
$\alpha\in A^{-n}S-A^{-n}S$, which gives $\alpha \in
[-2ra^{-n},2ra^{-n}]\times [-2rb^n,2rb^n]$. Let $\alpha=q\vu+p\vv$
for $p, q\in\Z$. Since $\vu,\vv$ are independent, there exists a
$C>0$ such that $|p|,|q| \leq Cb^n$.

Now write $\alpha=[x_0,y_0]^T$. Then $qu_1+pv_1=x_0$. Note that $|x_0| \leq 2ra^{-n}$.
It follows that
\begin{equation}  \label{3.2}
   \left|\frac{u_1}{v_1} +\frac{p}{q}\right| \leq \frac{2r|v_1|^{-1}}{|q|a^n}.
\end{equation}
Choose $n=m_k$ such that $Cb^{m_k} \leq M_k < Cb^{m_k+1}$ where $M_k$ are defined in
Lemma \ref{lem-3.3}. Write $a=b^{1+\ep}$, $\ep>0$. Then
$$
    a^{m_k}= \bigl(b^{m_k}\bigr)^{1+\ep} \geq (Cb)^{-1-\ep}M_k^{1+\ep}.
$$
Thus from (\ref{3.2}) we obtain
\begin{equation}  \label{3.3}
   \left|\frac{u_1}{v_1} +\frac{p}{q}\right| =
   \left|\frac{u_1}{v_1} -\frac{p_0}{|q|}\right| \leq
        \frac{2r|v_1|^{-1}(bC)^{1+\ep}}{|q|M_k^{1+\ep}},
\end{equation}
where $p_0/|q| = -p/q$. Note that $|q| \leq M_k$. By Lemma \ref{lem-3.3} we have
$$
   \left|\frac{u_1}{v_1} -\frac{p_0}{|q|}\right| >
           \frac{2r|v_1|^{-1}(bC)^{1+\ep}}{|q|M_k^{1+\ep}}
$$
for all sufficiently large $k$. This is a contradiction.
\eproof

\begin{prop}   \label{prop-3.5}
    Let $A=\diag\,(\lambda_1, \lambda_2)$ with $|\lambda_1|>1>|\lambda_2|$ and
$|\lambda_1\lambda_2|>1$. Let $\L=\Z\vu+\Z\vv$ be a full rank lattice in $\R^2$
with $\vu=[u_1,u_2]^T$ and $\vv=[v_1,v_2]^T$. Assume that $u_1/v_1\in\Q$.
Then there exists no measurable $T\subset\R^2$ such that
$T$ tiles translationally by $\L$ and multiplicatively by $A$.
\end{prop}
\Proof  We first consider the case $\vu=[1, c]^T$ and $\vv=[0,1]^T$. Then $F=(-1/2, 1/2]^2$
is a fundamental domain of $\L$, i.e. it tiles by $\L$. We also observe that the set
$\Omega =\{[x_1,x_2]^T:~|\lambda_1|^{-1}<2|x_1|\leq 1, x_2\in\R\}$ tiles multiplicatively by $A$.

Assume that a measurable $T\subset\R^2$ exists such that
$T$ tiles translationally by $\L$ and multiplicatively by $A$. Because $T$ tiles by $\L$ there
is a map $p_F:T\lra F$, which is one-to-one a.e.,
such that $x=p_F(x)+\alpha(x)$ for some unique $\alpha(x)\in\L$.
We have already introduced the address map $\tau_T$ induced by $T$ in (\ref{2.1}).
since both $\Omega$ and $T$ tile by $A$, $\tau_T:\Omega\lra T$ is bijective a.e..
Let $\phi: \Omega \lra F$ be given by $\phi=p_F\circ\tau_T$.

By definition of $\tau_T$ we have $\tau_T(x) = A^{n(x)}x$ for some $n(x)\in\Z$. Let
$$
     \Omega_n =\{x\in \Omega:~\tau_T(x) = A^n x\}.
$$
Then $\{\Omega_n:~n\in\Z\}$ is a partition of $\Omega$. Since translations are measure
preserving we have $\mu(\phi(\Omega_n)) =|\det(A)|^n\mu(\Omega_n)$.
It follows that for $n\geq 0$
\begin{equation}  \label{3.4}
    \sum_{n=0}^\infty \mu(\Omega_n) = \sum_{n=0}^\infty |\det(A)|^{-n}\mu(\phi(\Omega_n))
    \leq \sum_{n=0}^\infty \mu(\phi(\Omega_n))\leq 1.
\end{equation}
For $n<0$ observe that $A^n [x_1,x_2]^T= [\lambda_1^n x_1,
\lambda_2^n x_2]^T$. Now $|x_1|\leq 1/2$ for any $[x_1, x_2]^T\in
\Omega_n$. So $p_F(A^n [x_1,x_2]^T) =[\lambda_1^n x_1, z_2]^T$ for
some $|z_2| \leq 1/2$, and thus
$$
     \phi(\Omega_n) \subseteq [-|\lambda_1|^n/2, |\lambda_1|^n/2]\times [-1/2, 1/2].
$$
Hence $\mu(\phi(\Omega_n)) =|\det(A)|^n\mu(\Omega_n) \leq |\lambda_1|^n$, which yields
$\mu(\Omega_n) \leq |\lambda_2|^{-n}$. Thus
\begin{equation}  \label{3.5}
    \sum_{n<0}^\infty \mu(\Omega_n) \le  \sum_{n<0}^\infty |\lambda_2|^{-n} < \infty.
\end{equation}
Combining (\ref{3.4}) and (\ref{3.5}) yields $\mu(\Omega) < \infty$, a contradiction.

We now complete the proof of the proposition. In the general case, let $u_1/v_1\in\Q$.
Thus we can find a $\beta\in\R$ such that $\beta u_1=p$, $\beta v_1=q$ for some
$p,q\in\Z$, $\gcd(p,q)=1$. By Lemma ~\ref{lem-2.1} it suffices to prove there exists
no measurable $T$ such that $T$ tiles translationally by $\beta\L$ and multiplicatively
by $A$. Let $r,s\in\Z$ such that $rp+sq=1$, and let
$Q=\left[\begin{smallmatrix} r& -q \\s&p\end{smallmatrix}\right]$. Then
$\det(Q)=1$,
$$
    \beta\L =\begin{bmatrix} p & q\\\beta u_2 & \beta v_2\end{bmatrix}\Z^2
    =\begin{bmatrix} p & q\\\beta u_2 & \beta v_2\end{bmatrix}Q\Z^2
    =\begin{bmatrix} 1 & 0\\ b_1& b_2\end{bmatrix}\Z^2
$$
for some $b_1, b_2\in\R$. Finally, set $P_1=\diag(1, b_2^{-1})$. Then $P_1^{-1}AP_1=A$
and $P_1(\beta\L) = \Z[1, c]^T+\Z[0,1]^T$ where $c= b_1/b_2$, and it suffices to prove
there exists no measurable $T$ such that $T$ tiles translationally by
$\Z[1, c]^T+\Z[0,1]^T$ and multiplicatively by $A$. But this is precisely what
we have proved earlier.
\eproof

\vspace{1ex}
\noindent
{\bf Proof of Theorem \ref{theo-1.3}.}~~Let $\L=P\Z^2$ with $P\in\GLtwo$.
By Lemma~\ref{lem-2.1} there exists a measurable set $T$ that tiles by $\L$ and $A$
if and only if there exists a measurable set $T_1$ that tiles by $\Z^2$ and
$P^{-1}AP$.

\noindent
(a)~~If $|\det(A)|=1$ then there exists no $T$
such that $\mu(T) <\infty$ and $T$ tiles multiplicatively by $A$. Thus any $T$ that tiles
by $A$ cannot tile translationally by $\L$.

\noindent
(b)~~By Propositions \ref{prop-3.1} and \ref{prop-3.2},
there exists a measurable $T_1\subset\R^2$ such that $T_1$ tiles by $\Z^2$ and
$P^{-1}AP$. So there exists a measurable set $T\subset\R^2$ that tiles by $\L$ and $A$.

\noindent (c)~~Let $B=P^{-1}AP$. In this case both $\lambda_1$ and
$\lambda_2$ are real. Hence there exists a $Q\in\GLtwo$ such that
$Q^{-1}BQ = \diag(\lambda_1, \lambda_2)$. By Lemma \ref{lem-2.1}
the existence of a measurable $T_1$ that tiles $\R^2$
translationally by $\Z^2$ and multiplicatively by $B$ is
equivalent to the existence of a measurable $T_2$ such that $T_2$
tiles $\R^2$ translationally  by $Q^{-1}\Z^2$ and multiplicatively
by $Q^{-1}BQ=\diag (\lambda_1, \lambda_2)$. Let $Q=[\vu,\vv]$
where $\vu,\vv$ are the columns of $Q$, which are the eigenvectors
of $B=P^{-1}AP$ for the eigenvalues $\lambda_1$ and $\lambda_2$,
respectively. Hence
$$
    Q=\begin{bmatrix}u_1&v_1\\u_2&v_2\end{bmatrix}, \mhsp
    Q^{-1}=\frac{1}{\det(Q)}\begin{bmatrix}v_2&-v_1\\-u_2&u_1\end{bmatrix}.
$$
By Propositions \ref{prop-3.4} and \ref{prop-3.5}, a measurable $T_2$ that tiles translationally
by $Q^{-1}\Z^2$ and multiplicatively by $\diag(\lambda_1, \lambda_2)$ if and only if
$v_1/v_2\in\R\setminus\Q$. This proves the theorem.
\eproof

\vspace{1ex}
\noindent
{\bf Proof of Theorem \ref{theo-1.2}.}~~Since a wavelet set in $\R^2$ is a measurable set
$T$ that tiles translationally by $\Z^2$ and multiplicatively by $A^T$,
this theorem is a clearly corollary of Theorem \ref{theo-1.3}.
\eproof


\begin{thebibliography}{99999}
\bibitem{BeLe99}
J. J. Benedetto and M. T. Leon, {\it The construction of multiple
dyadic minimally supported frequency wavelets on $\R^d$}, in {\em
The Functional And Harmonic Analysis Of Wavelets And Frames}, Amer.
Math. Soc., Providence, RI, 1999.

\bibitem{DaLa98}
X. Dai and D. R. Larson, {\it Wandering vectors for unitary systems
and orthogonal wavelets}, {\em Mem. Amer. Math. Soc.} {\bf 134}, no.
640 (1998).

\bibitem{DLS97}
X. Dai, D. R. Larson, and  D. Speegle, {\it Wavelet sets in $\R^n$},
{\em J. Fourier Anal. Appl.} {\bf 3} (1997),  451--456

\bibitem{DLS97-b} X. Dai, D. Larson and D. Speegle,
{\it Wavelet in $\R^n$. II. }, {\em Wavelets, multiwavelets, and
their applications $($San Diego, CA, 1997$)$}, 15--40, Contemp.
Math., 216, Amer. Math. Soc., Providence, RI, 1998.

\bibitem{Dau88}
I.~Daubechies, {\it Orthonormal bases of compactly supported
wavelets}, {\em Comm. Pure Appl. Math.} {\bf 41} (1988), 906--996.

\bibitem{Dau-book}
I.~Daubechies,
{\em Ten Lectures on Wavelets}, SIAM, Philadelphia 1991.

\bibitem{FaWe96}
X. Fang and X.-H. Wang, {\it Construction of minimally supported
frequency wavelets}, {\em  J. Fourier Anal. Appl.} {\bf 2} (1996),
315--327.

\bibitem{HaWa01} D.-G. Han and Y. Wang,
{\it Lattice Tiling and the Weyl-Heisenberg families}, {\it
Geometric and Funct. Anal.} {\bf 11} (2001), no. 4, 742--758.

\bibitem{HWW96}
E. Hern{\'a}ndez, X.-H. Wang and G. Weiss, {\it Smoothing minimally
supported frequency wavelets. I}, {\em J. Fourier. Anal. Appl.} {\bf
2} (1996), 329--340.

\bibitem{HWW97}
E. Hern{\'a}ndez, X.-H. Wang and G. Weiss, {\it Smoothing minimally
supported frequency wavelets. II}, {\em J. Fourier. Anal. Appl.}
{\bf 3} (1997), 23--41.

\bibitem{HeWe-book}
E. Hern{\'a}ndez and G. Weiss, {\em A first course on wavelets},
Studies in Advanced Mathematics, CRC Press, Boca Raton, FL, 1996.

\bibitem{Ion02} E. Ionascu,
{\it A new construction of wavelet sets}, {\em Real Anal.
Exchange},{\bf 28} (2002/03), no. 2, 593--609.

\bibitem{ILP98}
E. J. Ionascu, D. R. Larson and C. M. Pearcy, {\it On wavelet sets},
{\em J. Fourier Anal. Appl.} {\bf 4} (1998), 711--721.

\bibitem{ILP98b} E. Ionascu, D. Larson, and  C. Pearcy,
{\it On the unitary systems affiliated with orthonormal wavelet
theory in n-dimensions}, {\em J. Funct. Anal.}, {\bf 157} (1998),
413--431.

\bibitem{IoPe98}
E. J. Ionascu and C. M. Pearcy, {\it On subwavelet sets}, {\em Proc.
Amer. Math. Soc.} {\bf 126} (1998), 3549--3552.

\bibitem{JaMa02} S. Jackson and R. D. Mauldin,
{\it On a lattice problem of H. Steinhaus}, {\em J. of AMS} {\bf 15}
(2002), 817--856.

\bibitem{Khi-book} A. Y. Khinchin, {\it Continued fractions}, Dover
Publications Inc., Mineola, New York, 1997.

\bibitem{Kol97} M. Kolountzakis, {\it Multi-lattice tiles}, {\em Internat. Math. Res. Notices} {\bf
19} (1997), 937--952.

\bibitem{LSST05} D. Larson, E. Schulz, D. Speegle and K. Taylor,
{\it Explicit cross sections of singly generated group actions}, to
appear.

\bibitem{LM86} P. G. Lemari$\acute{e}$ and Y. Meyer,
{\it Ondelettes et bases hilbertiennes}, {\em Rev. Math.
Iberoamericana}, {\bf 2}, (1986), 52--89.

\bibitem{M92} Y. Meyer,
{\it Wavelets and operators}, {\em Cambridge University Press},
(1992).

\bibitem{OlSp04} G. Olafsson and D. Speegle,
{\it Wavelets, wavelet sets, and linear actions on $\R^n$}, in {\em
Wavelets, Frames and Operator Theory}, 253-281, Contemp. Math. 345,
Amer. Math. Soc. Providence, RI, 2004.

\bibitem{SoWe98}
P. M. Soardi and D. Weiland, {\it Single wavelets in
$n$-dimensions}, {\em J. Fourier Anal. Appl.} {\bf 4} (1998),
299--315.

\bibitem{Spe03} D. Speegle,
{\it Dilation and translation tilings of $\R^n $ for non-expansive
dilations}, {\em Collect. Math.} {\bf 54} (2003), no. 2, 163--179.

\bibitem{Wan04} Y. Wang,
{\it Wavelets, tilings and spectral sets}, Duke Math J. {\bf 114}
(2004), no. 1, 43-57.
\end{thebibliography}
\end{document}